# The Mechanist's Challenge

## Did we really hope to get away with The Gödelian Argument?

Bhupinder Singh Anand[1]

## 1. Introduction – The Gödelian Argument

In a 1996 talk given at Oxford [Lu96], J.R. Lucas articulated his original, 1961, Gödelian argument [Lu61] afresh as follows:

> The Mechanist claims to have a model of the mind. We ask him whether it is consistent: if he cannot vouch for its consistency, it fails at the first examination; it just does not qualify as a plausible representation, since it does not distinguish those propositions it should affirm from those that it should deny, but is prepared to affirm both undiscriminatingly. We take the Mechanist seriously only if he will warrant that his purported model of the mind is consistent. In that case it passes the First Public Examination, but comes down at the Second, because knowing that it is consistent, we know that its Gödelian formula is true, which it cannot itself produce as true.

> More succinctly, we can, if a Mechanist presents us with a system that he claims is a model of the mind, ask him simply whether or not it can prove its Gödelian formula (according to some system of Gödel numbering). If he says it can, we know that it is inconsistent, and would be equally able to prove that 2 and 2 make 5, or that 0=1, and we waste little time on examining it. If, however, he acknowledges that the system cannot prove its Gödelian formula, then we know it is consistent, since it cannot prove every well-formed formula, and knowing that it is consistent, know also that its Gödelian formula is true.

[1] The author is an independent scholar. E-mail: re@alixcomsi.com; anandb@vsnl.com. Postal address: 32, Agarwal House, D Road, Churchgate, Mumbai - 400 020, INDIA. Tel: +91 (22) 2281 3353. Fax: +91 (22) 2209 5091.

In this formulation we have, essentially, a dialogue between the Mechanist and the Mentalist, as we may call him, with the Mechanist claiming to be able to produce a mechanist model of the Mentalist's mind, and the Mentalist being able to refute each particular instance offered.

## 2. Do we really hope to get away with this argument?

Do we really hope to get away with this argument?

For instance, when we ask a Mechanist whether or not his model can prove its Gödelian formula, the Mechanist would be quite justified in asking:

> Before I commit my resources to something that is impossible, can you assure me that it is provable?

Of course we cannot, since we, too, cannot give the proof of the formula that we are asking of the Mechanist's model.

If we were to, nevertheless, attempt to assert our superiority over the Mechanist by claiming, somewhat superciliously, that we do know, however, that the formula in question is true, but that the Mechanist's model does not, we would cut an even sorrier figure. The prompt, astonished - and subtly cunning - reply would be:

> Whatever gave you that idea! My model can even prove that it's true.

## 3. The Mechanist's Challenge

If, unwittingly, we were to express incredulity, we would face the humiliating challenge:

> The Gödelian formula ([Go31], p25, Gödel refers to it by its Gödel-number, 17*Genr*) is of the form [(A*x*)*R*(*x*)]. Neither of us can prove it from the agreed set of premises by the specified rules of inference of the Arithmetic in question ([Go31], p25(*1*), "17*Genr* is not *κ*-provable"). However, for a large enough numeral [*n*], my model can always provide a proof sequence for [*R*(1)&*R*(2)& ... &*R*(*n*)] from the premises faster than you.

The reason we would have to sheepishly concede defeat, without even addressing the challenge, is that Gödel ([1931], p22(*45*)) has given a primitive recursive relation, $xBy$[2], which ensures that the Mechanist - whose model contains Turing's Universal Machines as sub-models - can provide a proof of [*R*(*n*)] for any given numeral [*n*] (cf. [Go31], p26(*2*)).

Further, since our definition of the truth of a formula under an interpretation is based on Tarski's specification that it hold exhaustively, i.e., every instantiation of it must hold in the interpretation ([Me64], p51), our argument for our intellectual superiority over the Mechanist's model fails miserably in this particular instance.

End of story?

Not quite.

---

[2] Gödel ([Go1931], p22, def. 45) defines $xBy$ (primitive) recursively in the standard interpretation of a Peano Arithmetic so that it holds for natural numbers $x$, $y$, if, and only if, $x$ is the Gödel-number of a proof sequence in the Arithmetic for the formula of the Arithmetic whose Gödel-number is $y$.

## 4. There's more to come …

We could try to find solace, and console ourselves, with the fact that our brains are not designed ideally to tackle problems that can only be tackled by brute force, but are, surely, superior to the Mechanist's model when it comes to proofs that require ingenuity.

However, we will, then, face even greater discomfiture.

Gödel has given ([Go31], p22(*44*)) another primitive recursive predicate, $Bw(x)$[3], that the Mechanist can program on a UTM to - given sufficient time - create a set of theorems that will contain, and continue to be larger than, the set of all the arithmetical theorems that we have ever proven, and verified, collectively, using only pencil, paper, and our superior intellects!

Moreover, the set would contain theorems that we could not even conceptualise, since we would be unable to grasp even the statement, leave alone the proof, of these theorems.

## 5. Conclusion

So, what is really wrong with the Gödelian argument?

Well, ironically - and perhaps paradoxically - it's simply another instance of human insecurity applying very common, and very human, double standards - namely, matching our truth against the Mechanist's provability - to subtly load the dice in our favour when faced with semantic arguments that threaten to expose the continuing inability of the human ego to come to terms adequately with, and tap constructively into, the awesome conceptual power, diversity, and potential of any Intelligence, human or otherwise.

---

[3] Gödel ([Go1931], p22, def. 44) defines $Bw(x)$ (primitive) recursively in the standard interpretation of a Peano Arithmetic so that it holds for any natural number $x$ if, and only if, $x$ is the Gödel-number of a proof sequence in the Arithmetic.

(*Updated: Sunday*, *12<sup>th</sup> June 2005*, *4:19:22 PM IST*, *by re@alixcomsi.com*)